\documentclass[12pt]{amsart}

\usepackage{url}
\usepackage{hyperref}
\usepackage{amsmath,amssymb,latexsym, amsfonts, amscd, amsthm, 
xy,verbatim}

\xyoption{all}

\theoremstyle{plain}

\newtheorem{theorem}{Theorem}[section]
\newtheorem{conjecture}[theorem]{Conjecture}

\theoremstyle{definition}

\newtheorem{remark}[theorem]{Remark}

\long\def\symbolfootnote[#1]#2{\begingroup
\def\thefootnote{\fnsymbol{footnote}}\footnote[#1]{#2}\endgroup}

\def\SL{{\bf SL}}

\def\PSL{{\bf PSL}}

\def\1{{\mathbf 1}}

\DeclareMathOperator{\MM}{M}

\def\bdf{\begin{defn}}
\def\edf{\end{defn}}

\def\cO{\mathcal{O}}

\newcommand{\pty}[1]{\langle #1 \rangle}
\newcommand{\RR}{\mathbf{R}}
\newcommand{\ZZ}{\mathbf{Z}}
\newcommand{\CC}{\mathbf{C}}
\newcommand{\QQ}{\mathbf{Q}}

\newcommand{\smallabcd}{\left(\begin{smallmatrix}a&b\\c&d\end{smallmatrix}\right)}
\newcommand{\abcd}{\begin{pmatrix}a&b\\c&d\end{pmatrix}}

\renewcommand{\epsilon}{\varepsilon}
\renewcommand{\phi}{\varphi}

\begin{document}
\title{Rationality of secant zeta values}

\author{Pierre Charollois}

\email{charollois@math.jussieu.fr}
\address{Institut de math\'ematiques de Jussieu, Universit\'e  Paris 6, 4 Place Jussieu, 75005 Paris, France}
\author{Matthew Greenberg}
\email{mgreenbe@ucalgary.ca}
\address{Department of Mathematics and Statistics, University of Calgary, 2500 University Drive NW, Calgary, Alberta, T2N 1N4, Canada}

\date{\today}
\thanks{PC's research is partially supported by the grant REGULATEURS  ``ANR-12-BS01-0002''.} 
\thanks{MG's research is supported by NSERC of Canada}

\begin{abstract} We use the Arakawa-Berndt theory of generalized $\eta$-functions to prove a conjecture of Lal\'in, Rodrigue and Rogers concerning the algebraic nature of special values of the secant zeta  function.
\end{abstract}
\maketitle

\section{Introduction}

The cotangent and secant zeta functions attached to an  algebraic irrational  number $\alpha$ are defined, for $\textrm{Re}(s)$ large enough, by the Dirichlet series 
\[
\xi(\alpha,s)=\sum_{n=1}^\infty \frac{\cot(\pi n \alpha)}{n^s}\qquad \text{and}\qquad \psi(\alpha,s)=\sum_{n=1}^\infty\frac{\sec(\pi n \alpha)}{n^s},
\]
respectively. The cotangent zeta function was introduced for $s=2k+1$, $k\geq 1$ an integer, by Lerch~\cite{Ler}, and  for general $s$ by Berndt~\cite{Ber76} in the course of his study of generalized Dedekind sums. 
Lerch stated (without proof) the following functional equation valid for algebraic irrational $\alpha$ and sufficiently large $k=k(\alpha)$:
\begin{equation}\label{E:LerchFE}
\xi(\alpha,2k+1) + \alpha^{2k} \xi(\tfrac1\alpha,2k+1)=(2\pi)^{2k+1}\phi(\alpha,2k+1),
\end{equation}
where
\[
\phi(\alpha,n)=\sum_{j=0}^{n+1}\frac{B_jB_{n+1-j}}{j!(n+1-j)!}\alpha^{j-1},
\]
and  $B_k$ is the $k$-th  Bernoulli  number.
Suppose $\alpha\neq \pm 1$ is a unit in the quadratic field $\QQ(\sqrt{d})$ of discriminant $d$.  Writing $\alpha=\tfrac{a+b\sqrt{d}}{2}$ and defining $\epsilon=\pm 1$ by $a^2-b^2d=\epsilon$, we have $\tfrac1\alpha=-\epsilon(\alpha-a)$.  Using the obvious identities $\xi(-\alpha,s)=-\xi(\alpha,s)$ and $\xi(\alpha-1,s)=\xi(\alpha,s)$ together with~\eqref{E:LerchFE}, we deduce the following rationality result:
\begin{equation}\label{E:LerchRat}
\frac{\xi(\alpha,2k+1)}{(2\pi)^{2k+1}\sqrt{d}}=\frac{1}{\sqrt{d}}\frac{\phi(\alpha,2k+1)}{1-\epsilon\alpha^{2k}}\in\QQ. 
\end{equation}
Another proof of this result was given by Berndt~\cite[Theorem 5.2]{Ber76}. 
More generally, if $\alpha$ is an arbitrary real quadratic irrationality, then Lerch uses the continued fraction expansion of $\alpha$ to conclude that the left hand side of~\eqref{E:LerchRat} belongs to $\QQ(\alpha)$ -- see~\cite[(3)]{Ler}.



Lal\'in, Rodrigue and Rogers conjecture an analogous result for the secant zeta function:
\begin{conjecture}[{\cite[Conjecture 1]{LRR}}]\label{C:LRR}
Suppose $d$ and $k$ are positive integers such that $d$ is not a square.  Then
\[
\frac{\psi(\sqrt{d},2k)}{\pi^{2k}}\in\QQ.
\]
\end{conjecture}
\noindent They prove a functional equation for~$\psi(\alpha,s)$ analogous to~\eqref{E:LerchFE} and use it to deduce many instances of Conjecture~\ref{C:LRR} as above.  In \S\ref{S:LRRFE}, we show that their functional equation is actually sufficient to prove in general that $\psi(\alpha,2k)\in\pi^{2k}\QQ(\alpha)$ for all real quadratic irrationalities $\alpha$, not merely those of the form $\alpha=\sqrt{d}$. 

 In \S\ref{S:eta}, we relate $\psi(\alpha,s)$ to the \emph{generalized $\eta$-functions} studied by Arakawa~\cite{Ar82}, fascinating objects in their own right.  Using this relationship, we leverage Arakawa's results to give an explicit formula for $\psi(\alpha,2k)$ for real quadratic irrationalities $\alpha$, yielding another proof of Conjecture~\ref{C:LRR}.

\section{The Lal\'in-Rodrigue-Rogers functional equation}\label{S:LRRFE}
One can prove Conjecture~\ref{C:LRR} following Lerch's approach for the cotangent zeta function.  Most of this argument was given in~\cite{LRR}; we complete their thought.  Let $A,B\in\PSL_2(\ZZ)$ be defined by
\[
A=\begin{pmatrix}1&2\\0&1\end{pmatrix},\qquad B=\begin{pmatrix}1&0\\2&1\end{pmatrix}.
\]
The following functional equations are established in~\cite[(4.1), (4.2)]{LRR}:
\begin{align*}
\psi(A\alpha,2k)&=\psi(\alpha,2k),\\
\psi(B\alpha,2k)&=(2\alpha+1)^{1-2k}\psi(\alpha,2k)  \\
&\quad - \frac{\pi^{2k}}{(2k)!}\sum_{m=0}^{2k}(2^{m-1}-1)B_mE_{2k-m}\binom{2k}{m}(\alpha+1)^{2k-m}\Big((2\alpha+1)^{m-2k} - (2\alpha+1)^{1-2k}\Big).
\end{align*}
Here, $B_n$ is the $n$-th Bernoulli number and $E_n$ is the $n$-th Euler number.  It follows that if $C\in\langle A,B\rangle$ then there is a $\QQ(\alpha)$-linear relation between $\psi(C\alpha,2k)$, $\psi(\alpha,2k)$ and $\pi^{2k}$.  Thus, if there is a matrix $C\in\langle A,B\rangle$ such that $C\alpha=\alpha$ then $\psi(\alpha,2k)\in\pi^{2k}\QQ(\alpha)$.  In~\cite[\S4]{LRR}, several families of examples of such pairs $(\alpha,C)$ are given and the associated linear relations are worked out explicitly.  We merely point out that if $\alpha$ is any real quadratic irrationality, then there is \emph{always} a $C\in\langle A,B\rangle$, $C\neq 1$ such that $C\alpha=\alpha$.

To see this, let $\alpha$ be a real quadratic irrationality and consider the lattice $L=\ZZ+\ZZ\alpha\subset\QQ(\alpha)$.  Let $\cO$ be the order of $L$:
\[
\cO=\{x\in \QQ(\alpha) : xL\subset L\}.
\]
Then $\cO$ is an order in $\QQ(\alpha)$.  Let  $u\in \cO.$
Writing $u\cdot\alpha=a\alpha+b$ and $u\cdot 1=c\alpha+d$ with $a,b,c,d\in\ZZ$, we have
\[
u\begin{pmatrix}\alpha\\1\end{pmatrix}=\abcd\begin{pmatrix}\alpha\\1\end{pmatrix}.
\]
Write $j(u)$ for the matrix $\smallabcd$ associated to $u$ above.  Then $j:\cO\to\MM_2(\ZZ)$ is a ring homomorphism.  Since $\left(\begin{smallmatrix}\alpha\\1\end{smallmatrix}\right)$ is an eigenvector of $j(u)$ for all $u\in\cO$, we have $j(u)\alpha:=\frac{a\alpha+b}{c\alpha+d}=\alpha$ when $u\neq 0$.
By Dirichlet's unit theorem, the group $\cO_+^*$ of totally positive units in $\cO$ is free of rank $1$; write $\cO^*_+=\langle\gamma\rangle$.
By \cite[p.\ 84]{Sc74},  $\langle A,B\rangle$ is the principal congruence subgroup $\Gamma(2)\subset\PSL_2(\ZZ)$, this inclusion having index 6.  Therefore, $C:= j(\gamma^6)$ satisfies $C\neq 1$, $C\in\Gamma(2)$ and $C\alpha=\alpha$.

\section{Generalized $\eta$-functions and secant zeta values}\label{S:eta}

Arakawa~\cite{Ar88} gave another proof of  (\ref{E:LerchRat}) by relating $\xi(\alpha,s)$ to generalized $\eta$-functions, the theory of which he developed in~\cite{Ar82}.  In turn, Arakawa's work has its foundations in papers of Lewittes~\cite{Lew72} and Berndt~\cite{Ber73}.  We show that Arakawa's method can also be used to analyze the secant zeta function.

For $x\in\RR$, define $\pty{x}$ (resp., $\{x\}$) by
\[
0<\pty{x}\leq 1\qquad \text{(resp.\ $0\leq \{x\}< 1)$}\qquad  and \qquad x-\pty{x}\in\ZZ\qquad\text{(resp., $x-\{x\}\in\ZZ$)}.
\]
We set $e(z)=e^{2\pi i z}$.

Following~\cite{Ar82}, let $p,q\in\RR$ and define
\begin{align*}
\eta(\alpha,s,p,q)&=\sum_{n=1}^\infty n^{s-1}\frac{e(n(p\alpha+q))}{1-e(n\alpha)}\\
H(\alpha,s,p,q) &= \eta(\alpha,s,\pty{p},q) + e\left(\tfrac{s}{2}\right)\eta(\alpha,s,\pty{-p},-q).
\end{align*}

\begin{theorem}[{\cite[Lemma 1 and Theorem 2]{Ar82}}] Suppose $$\alpha\in\RR\cap\bar{\QQ}\qquad\text{and}\qquad \alpha\notin\QQ.$$  Then $\eta(\alpha,s,p,q)$ is absolutely convergent for $\Re(s)<0$.  If, in addition, $$[\QQ(\alpha):\QQ]=2\qquad\text{and}\qquad p,q\in\QQ$$ then $H(\alpha,s,p,q)$ has analytic continuation to $\CC-\{0\}$, and the singularity at $s=0$ is at worst a simple pole.
\end{theorem}

\begin{remark}
The convergence of $\eta(\alpha,s,p,q)$ relies on the Thue-Siegel-Roth theorem in much the same way that the convergence of $\psi(\alpha,s)$ does  -- see~\cite[Theorem 1]{LRR}.
\end{remark}

An elementary computation yields
\[
\xi(\alpha,s)=-2i\left(\frac{H(\alpha,1-s,1,0)}{1+e(\frac s 2)}+\frac{1}{2}\zeta(s)\right),
\]
so rationality statements for $H(\alpha,1-s,1, 0)$ for even integral $s$ translate to rationality statements for $\xi(\alpha,s)$ at odd integral $s$.  In contrast, $\psi(\alpha,s)$ does not seem to have a simple expression in terms of $H(\alpha,s,p,q)$.  Crucially, however, we still have a relation between certain special values of $H$ and $\psi$ relying on the relation $\frac 34=-\frac 14+1$ :
\begin{align*}
\psi\left(\tfrac \alpha 2,1-s\right)&=\sum_{n=1}^\infty n^{s-1}\frac{2}{e(\frac{n\alpha}4)+e(-\frac{n\alpha}4)}\\
&=2\sum_{n=1}^\infty n^{s-1}\frac{e(\frac {n\alpha}4)}{1-e(n\alpha)} - 2\sum_{n=1}^\infty n^{s-1}\frac{e(\frac{3n\alpha}4)}{1-e(n\alpha)}\\
& =2\eta(\alpha,s,\pty{\tfrac 14},0)-2\eta(\alpha,s,\pty{-\tfrac 14},0).
\end{align*}

\noindent If $s=1-2k$, then $e(\frac s 2)=-1$ and we conclude that 
\begin{equation}\label{E:psiH}
\psi(\tfrac \alpha 2,2k)=2H(\alpha,1-2k,\tfrac14,0).
\end{equation}

For the rest of the paper, suppose that $\alpha$ is a real quadratic irrationality. Formulas of Berndt and Arakawa allow us to evaluate $H(\alpha,1-2k,\frac 14,0)$ rather explicitly.  Let \[V=\abcd\in\SL_2(\ZZ)\] be a 
matrix such that 
\[
c>0\qquad \text{and}\qquad \beta:=c\alpha+d>0.
\]
Set
\[
(p',q')=(p,q)V\qquad\text{and}\qquad \varrho=\{q'\}c-\{p'\}d.
\]

\begin{theorem}[{\cite[Theorem 1 and Eq. (1.19)]{Ar82}}]\label{T:trans}
Suppose that $p$ and $p'$ are not in $\ZZ$. Then the following transformation formula holds:
\[
\beta^{-s}H(V\alpha,s,p,q) = H(\alpha,s,p',q')+(2\pi)^{-s}e(-\frac s 4)L(\alpha,s,p',q',c,d),
\]
where $L(\alpha,s,p',q',c,d)$ is as in (2) of~\cite{Ber73}.
If $s=-m$ is a negative integer, then
\[
L(\alpha,-m,p',q',c,d)=2\pi i\sum_{j=1}^c\sum_{\ell=0}^{m+2}
b_\ell\left(\frac{j-\{p'\}}{c}\right)b_{m+2-\ell}\left(\left\{\frac{jd+\varrho}{c}\right\}\right)(-\beta)^{\ell-1}.
\]
Here, $b_\ell$ is the (normalized)  $\ell$-th Bernoulli polynomial defined by the generating series $$\frac{ue^{ux}}{e^u-1}=\sum_{\ell\geq 0} b_\ell(x) u^\ell.$$
\end{theorem}

By~\cite[Lemma 4]{Ar82}, for any rational numbers $ p$ and $q$ there is a totally positive unit $\beta$ of $\QQ(\alpha)$ and a matrix $V\in\SL_2(\ZZ)$ such that 
\[
c>0,\qquad (p',q'):=(p,q)V\equiv (p,q)\pmod{\ZZ^2},\qquad \text{and}\qquad V\begin{pmatrix}\alpha\\1\end{pmatrix}=\beta\begin{pmatrix}\alpha\\1\end{pmatrix}. 
\]
The last condition implies that $\beta=c\alpha+d$, consistently with the notation introduced above.  
Suppose $p\notin\ZZ$ .  Then $p'\notin\ZZ$, too, as $p\equiv p'\pmod{\ZZ}$.
Applying Theorem~\ref{T:trans}, observing that $H$ and $L$  only depend on the class of $(p,q)\equiv (p',q')$ modulo $\ZZ^2$, and rearranging terms, we get
\[
(\beta^{-s}-1)H(\alpha,s,p,q)=(2\pi)^{-s}e(-\tfrac s4)L(\alpha,s,p,q,c,d).
\]
By the second part of Theorem~\ref{T:trans}, if $s=1-2k$ then
\begin{multline}
\frac{H(\alpha,1-2k,p,q)}{\pi^{2k}}=\frac{2^{2k}(-1)^{k}}{(\beta^{2k-1}-1)}\sum_{j=1}^c\sum_{\ell=0}^{2k+1}
b_\ell\left(\frac{j-\{p\}}{c}\right)b_{2k+1-\ell}\left(\left\{\frac{jd+\varrho}{c}\right\}\right)(-\beta)^{\ell-1}.
\end{multline}

\noindent Setting $(p,q)=(\frac14,0)$  and using~\eqref{E:psiH},  this formula specializes  to 

\begin{equation}
\frac{\psi(\frac \alpha 2,2k)}{\pi^{2k}}=\\ \frac{2^{2k+1}(-1)^{k} }{(\beta^{2k-1}-1)}\sum_{j=1}^c\sum_{\ell=0}^{2k+1}
b_\ell\left(\frac{j-\frac 14}{c}\right)b_{2k+1-\ell}\left(\left\{\frac{d(j-\frac 14)}{c}\right\}\right)(-\beta)^{\ell-1}.
\end{equation}
We conclude :
\begin{theorem}
\label{T:main}
Suppose $\alpha$ is a real quadratic irrationality and $k$ is a positive integer. Then
\begin{equation*}
\frac{\psi(\alpha,2k)}{\pi^{2k}}\in\QQ(\alpha).
\end{equation*}
Moreover, if $x\mapsto x'$ is the nontrivial automorphism of $\QQ(\alpha)$, then
\begin{equation*}
\left(\frac{\psi(\alpha,2k)}{\pi^{2k}}\right)'=\frac{\psi(\alpha',2k)}{\pi^{2k}}.
\end{equation*}
\end{theorem}

\noindent Conjecture~\ref{C:LRR} follows from Theorem~\ref{T:main} and the evenness of the secant function.

\end{document}